\documentclass[oneside,english]{amsart}
\usepackage[T1]{fontenc}
\usepackage[latin9]{inputenc}
\usepackage{float}
\usepackage{calc}
\usepackage{amsthm}
\usepackage{graphicx}
\usepackage{esint}
\usepackage[pdftex,bookmarks=true,%
bookmarksdepth=3,%
bookmarksnumbered=true,%
bookmarkstype=toc]{hyperref}
\makeatletter

\numberwithin{equation}{section}
\numberwithin{figure}{section}
\theoremstyle{plain}
\newtheorem{thm}{\protect\theoremname}
  \theoremstyle{plain}
  \newtheorem{cor}[thm]{\protect\corollaryname}
  \theoremstyle{remark}
  \newtheorem{rem}[thm]{\protect\remarkname}

\makeatother

\usepackage{babel}
  \providecommand{\corollaryname}{Corollary}
  \providecommand{\remarkname}{Remark}
\providecommand{\theoremname}{Theorem}

\begin{document}

\title{The Zetafast algorithm for computing zeta functions}

\author{Kurt Fischer}

\address{Tokuyama College of Technology, Gakuendai, Shunan, Yamaguchi 745-8585,
JAPAN}


\begin{abstract}
We express the Riemann zeta function $\zeta\left(s\right)$ of argument
$s=\sigma+i\tau$ with imaginary part $\tau$ in terms of three absolutely
convergent series. The resulting simple algorithm allows to compute,
to arbitrary precision, $\zeta\left(s\right)$ and its derivatives
using at most $C\left(\epsilon\right)\left|\tau\right|^{\frac{1}{2}+\epsilon}$
summands for any $\epsilon>0$, with explicit error bounds. It can
be regarded as a quantitative version of the approximate functional
equation. The numerical implementation is straightforward. The approach
works for any type of zeta function with a similar functional equation
such as Dirichlet $L$-functions, or the Davenport-Heilbronn type
zeta functions. 
\end{abstract}

\maketitle

\section{Results\label{sec:Results}}
\begin{thm}
Zetafast algorithm for Riemann zeta function:\label{thm:Zetafast algorithm for Riemann zeta function}

Given a positive integer $v$, $\zeta\left(s\right)$ has for $\sigma\leq\sigma_{0}<v$
and $s\notin\left\{ 1,\dots,v-1\right\} $ the following representation
in terms of three absolutely and uniformly converging series $D\left(s\right)$,
$E_{1}\left(s\right)$, and $E_{-1}\left(s\right)$,

\begin{eqnarray}
\zeta\left(s\right) & = & D\left(s\right)+\sum_{\mu=\pm1}E_{\mu}\left(s\right)-\frac{\Gamma\left(1-s+v\right)}{\left(1-s\right)\Gamma\left(v\right)}N^{1-s}\label{eq:Zetafast in terms of D(S) and E_mu(s)}\\
D\left(s\right) & = & \sum_{n=1}^{\infty}n^{-s}Q\left(v,\frac{n}{N}\right)\label{eq:D(s)}\\
E_{\mu}\left(s\right) & = & \left(2\pi\right)^{s-1}\Gamma\left(1-s\right)e^{i\mu\frac{\pi}{2}\left(1-s\right)}\sum_{m=1}^{\infty}E_{\mu}\left(m,s\right)\label{eq:E_mu(s)}\\
E_{\mu}\left(m,s\right) & = & m^{s-1}-\sum_{w=0}^{v-1}\binom{s-1}{w}\left(m+\frac{i\mu}{2\pi N}\right)^{s-1-w}\left(\frac{-i\mu}{2\pi N}\right)^{w}\label{eq:E_mu(m,s)}
\end{eqnarray}
where $Q\left(v,m\right)$ is the normalized incomplete gamma function
\begin{align}
Q\left(v,m\right) & =\sum_{w=0}^{v-1}\frac{m^{w}}{w!}e^{-m}\label{eq:Normalized incomplete gamma function}
\end{align}
For positive integer $s=k$ with $k<v$ one has to take the limit
$\lim_{s\to k}E_{\mu}\left(s\right)$.\end{thm}
\begin{cor}
The derivatives of $\zeta\left(s\right)$ can be calculated by differentiating
the above series termwise. 
\end{cor}
We are mainly interested in the critical strip and its surroundings,
for which we have the following rough but explicit estimate:
\begin{thm}
Accuracy and speed of Zetafast algorithm:\label{thm:Accuracy and speed of Zetafast algorithm}

For a given argument $s$ with 
\begin{eqnarray}
0 & \leq\sigma\leq & 2\label{eq:0<=00003D sigma <=00003D 2}\\
0 & < & \tau\label{eq:tau >=00003D0}
\end{eqnarray}
and accuracy 
\begin{eqnarray}
\delta & \leq & 0.05\label{eq:delta max}
\end{eqnarray}
choose $v$ as the next higher integer $\left\lceil x_{0}\right\rceil $
of the unique solution $x_{0}$ of
\begin{eqnarray}
x-\max\left(\frac{1-\sigma}{2},0\right)\ln\left(\frac{1}{2}+x+\tau\right) & = & \ln\frac{8}{\delta}\label{eq:v approx best condition log version again}
\end{eqnarray}
in the unknown $x$, and from this $M=\left\lceil N\right\rceil $
with 

\begin{eqnarray}
N & = & 1.11\left(1+\frac{\frac{1}{2}+\tau}{v}\right)^{\frac{1}{2}}\label{eq:N sufficient}
\end{eqnarray}
as well as
\begin{eqnarray}
\lambda & = & 3.151\label{eq:lambda sufficient}
\end{eqnarray}
Then we have for some real number $c$ with $\left|c\right|<1$ the
approximation 
\begin{eqnarray}
\zeta\left(s\right) & = & D\left(N,s\right)+E_{1}\left(M,s\right)-\frac{\Gamma\left(1-s+v\right)}{\left(1-s\right)\Gamma\left(v\right)}N^{1-s}+c\delta\label{eq:Zetafast approximation}\\
D\left(N,s\right) & = & \sum_{n=1}^{\left\lceil \lambda vN\right\rceil }n^{-s}Q\left(v,\frac{n}{N}\right)\nonumber \\
E_{1}\left(M,s\right) & = & \left(2\pi\right)^{s-1}\Gamma\left(1-s\right)e^{i\frac{\pi}{2}\left(1-s\right)}\sum_{m=1}^{M}E_{1}\left(m,s\right)\nonumber 
\end{eqnarray}
Under the condition 
\begin{eqnarray}
\tau & > & \frac{5}{3}\left(\frac{3}{2}+\ln\frac{8}{\delta}\right)\label{eq:tau large condition}
\end{eqnarray}
we need at most $S$ summands to calculate $D\left(N,s\right)$ and
$E_{1}\left(M,s\right)$, with

\begin{eqnarray}
S & = & 2+8\sqrt{1+\ln\frac{8}{\delta}+\max\left(\frac{1-\sigma}{2},0\right)\ln\left(2\tau\right)}~\sqrt{\tau}\label{eq:Zetafast speed}
\end{eqnarray}
For $\tau=0$, the same estimates hold, but we can neglect $E_{1}\left(M,s\right)$
altogether.\end{thm}
\begin{rem}
Estimate~(\ref{eq:Zetafast speed}) shows that Zetafast allows for
arbitrary precision, while being essentially as fast as the Riemann-Siegel
formula~\cite{Titchmarsh_Zeta}. Using similar arguments, we can
obtain explicit error bounds for the derivatives of the Riemann zeta
function, or for Dirichlet $L$-functions or its derivatives. What
is more, there is room for further tightening the error bounds, or
to accelerate the algorithm along the lines of~\cite{Hiary_Fast.methods.to.compute.Riemann.zeta},
or~\cite{Odlyzko_Schoenhage} for multiple evaluations.
\end{rem}
The same arguments show that linear combinations of Dirichlet $L$-functions
$L\left(\chi_{k},s\right)$ for $k=1,2,\dots,l$ such as the Davenport-Heilbronn
zeta function or Hurwitz zeta functions $\zeta\left(s,r\right)$ for
rational parameters $r$, have a Zetafast algorithm, if they obey
a functional equation expressing this linear combination in terms
of a linear combination of functions $r_{k}^{\left(1-s\right)}e^{i\phi_{k}\left(1-s\right)}\Gamma\left(1-s\right)L\left(\chi_{k},1-s\right)$
for real numbers $r_{k}$ and $\phi_{k}$ with $\left|\phi_{k}\right|<\pi$.
We give here only one example:
\begin{thm}
Zetafast algorithm for Dirichlet $L$-functions:\label{thm:Zetafast algorithm for Dirichlet L-functions}

A Dirichlet $L$-function with primitive, non-principal character
$\chi$ and Gauss sum $G\left(\chi\right)$ can be calculated for
$\sigma\leq\sigma_{0}<v$ and $s\notin\left\{ 1,\dots,v-1\right\} $
by the following absolutely converging series, 
\begin{align}
L\left(s,\chi\right) & =D\left(s,\chi\right)+\left(2\pi\right)^{s-1}\Gamma\left(1-s\right)q^{-s}G\left(\chi\right)\sum_{\mu=\pm1}\chi\left(-\mu\right)e^{i\mu\frac{\pi}{2}\left(1-s\right)}E_{\mu}\left(s,\chi\right)\label{eq:Zetafast Dirichlet L-function}\\
D\left(s,\chi\right) & =\sum_{n=1}^{\infty}\chi\left(n\right)n^{-s}Q\left(v,\frac{n}{N}\right)\nonumber \\
E_{\mu}\left(s,\chi\right) & =\sum_{m=1}^{\infty}\overline{\chi}\left(m\right)\left[m^{s-1}-\sum_{w=0}^{v-1}\binom{s-1}{w}\left(m+\frac{i\mu q}{2\pi N}\right)^{s-1-w}\left(\frac{-i\mu q}{2\pi N}\right)^{w}\right]\nonumber 
\end{align}
For positive integer $s=k$ with $k<v$ one has to take the limit
$\lim_{s\to k}E_{\mu}\left(s,\chi\right)$.
\end{thm}

\section{Proofs\label{sec:Proofs}}

\subsection{Proof of Theorem~\ref{thm:Zetafast algorithm for Riemann zeta function}\label{sub:Proof of Theorem=0000A0=00005Bthm:Zetafast=00005D}}

We assume for the moment $0<\sigma<1$, start with the series~(\ref{eq:D(s)})
and express the cutoff $Q$ in terms of its inverse Mellin transform, 

\begin{eqnarray}
D\left(s\right) & = & \sum_{w=0}^{v-1}\frac{1}{w!}\sum_{n=1}^{\infty}n^{-s}\int_{x>1-\sigma}\Gamma\left(z+w\right)\left(\frac{n}{N}\right)^{-z}\frac{dz}{2\pi i}\label{eq:D(w,s)}
\end{eqnarray}
where $x>1-\sigma$ denotes that the integration contour is for this
value of $x$ along the vertical line from $x-i\infty$ to $x+i\infty$.
We interchange the absolutely convergent series and integrals,
\begin{eqnarray*}
D\left(s\right) & = & \sum_{w=0}^{v-1}\frac{1}{w!}\int_{x>1-\sigma}\Gamma\left(z+w\right)N^{z}\zeta\left(z+s\right)\frac{dz}{2\pi i}
\end{eqnarray*}
Because the $\Gamma$ function decreases exponentially with increasing
$\left|y\right|$ while the $\zeta$ function increases at most algebraically~\cite{Titchmarsh_Zeta},
we can move the contour to $-1<x<-\sigma$ and pick up the residues
at $z+s=1$ and, for the case $w=0$, at $z=0$,
\begin{eqnarray}
D\left(s\right) & = & \sum_{w=0}^{v-1}\frac{\Gamma\left(1-s+w\right)}{w!}N^{1-s}+\zeta\left(s\right)+\sum_{w=0}^{v-1}\frac{1}{w!}E\left(w,s\right)\label{eq:D(w,s) with integral rest E(w,s)}\\
E\left(w,s\right) & = & \int_{-1<x<-\sigma}\Gamma\left(z+w\right)N^{z}\zeta\left(z+s\right)\frac{dz}{2\pi i}\label{eq:E(w,s)}
\end{eqnarray}
The idea is now to use the functional equation for $\zeta$, 
\begin{eqnarray}
\zeta\left(s+z\right) & = & \left(2\pi\right)^{s+z-1}\Gamma\left(1-s-z\right)\zeta\left(1-s-z\right)\sum_{\mu=\pm1}e^{i\mu\frac{\pi}{2}\left(1-s-z\right)}\label{eq:zeta functional equation}
\end{eqnarray}
and to re-express the $\zeta$ function in terms of its then absolutely
convergent Dirichlet series, because $1-\sigma-x>1$, so that

\begin{eqnarray}
E\left(w,s\right) & = & \left(2\pi\right)^{s-1}\Gamma\left(1-s\right)\sum_{\mu=\pm1}e^{i\mu\frac{\pi}{2}\left(1-s\right)}\sum_{m=1}^{\infty}E_{\mu}\left(w,m,s\right)\label{eq:E_mu in terms of Mellin-Barnes}\\
E_{\mu}\left(w,m,s\right) & = & m^{s-1}\int_{-1<x<-\sigma}\frac{\Gamma\left(z+w\right)\Gamma\left(1-s-z\right)}{\Gamma\left(1-s\right)}\left(\frac{i\mu}{2\pi Nm}\right)^{-z}\frac{dz}{2\pi i}\nonumber 
\end{eqnarray}
This Mellin-Barnes integral is the well-known inverse of the beta
function integral~(\cite{NIST_online}~5.12.3 and~5.13.1), thus 

\begin{align}
E_{\mu}\left(w,m,s\right) & =\frac{\Gamma\left(1-s+w\right)}{\Gamma\left(1-s\right)}\left(m+\frac{i\mu}{2\pi N}\right)^{s-1-w}\left(\frac{i\mu}{2\pi N}\right)^{w}-\delta_{w0}m^{s-1}\label{eq:Mellin-Barnes integral evaluated}
\end{align}
From~(\ref{eq:D(w,s) with integral rest E(w,s)}),~(\ref{eq:E(w,s)}),~(\ref{eq:E_mu in terms of Mellin-Barnes})
and~(\ref{eq:Mellin-Barnes integral evaluated}) and the elementary
identities, 
\begin{eqnarray}
\sum_{w=0}^{v-1}\frac{\Gamma\left(1-s+w\right)}{w!} & = & \frac{\Gamma\left(1-s+v\right)}{\left(1-s\right)\Gamma\left(v\right)}\label{eq:Binomial addition theorem}\\
\frac{\Gamma\left(1-s+w\right)}{\Gamma\left(1-s\right)w!} & = & \binom{s-1}{w}\left(-1\right)^{w}\label{eq:Binomial coefficient in terms of gamma functions}
\end{eqnarray}
and setting 
\begin{eqnarray*}
-\sum_{w=0}^{v-1}\frac{1}{w!}E_{\mu}\left(w,m,s\right) & = & E_{\mu}\left(m,s\right)
\end{eqnarray*}
we arrive at our representation~(\ref{eq:Zetafast in terms of D(S) and E_mu(s)})-(\ref{eq:E_mu(m,s)}). 

Finally, absolute and uniform convergence is obvious for $D\left(s\right)$,
and follows for $E_{\mu}\left(s\right)$ in the region $\sigma\le\sigma_{0}<v$
by analytic continuation because by the binomial theorem $\Gamma\left(1-s\right)E_{\mu}\left(m,s\right)\sim\Gamma\left(1-s+v\right)m^{s-1-v}$
as $m\to\infty$. This shows also that for $s=k$ for a positive integer
$k<v$, the poles of $\Gamma\left(1-s+w\right)$ for $w=0,1,\dots,s-1$
cancel each other, so that we can take the limit $\lim_{s\to k}E_{\mu}\left(s\right)$.

\subsection{Proof of Theorem~\ref{thm:Accuracy and speed of Zetafast algorithm}\label{sub:Proof of Theorem=0000A0=00005Bthm:Zetafast speed=00005D}}

\subsubsection{Choosing $v$\label{sub:Choosing v}}

Because we assume~(\ref{eq:0<=00003D sigma <=00003D 2}) and $\sigma<v$,
we are free to restrict $v$ to 
\begin{eqnarray}
v & \geq & 5\label{eq:v>=00003D5}
\end{eqnarray}
First we prove that equation~(\ref{eq:v approx best condition log version again})
has a unique solution $x_{0}>5$. 

For the case $\sigma\geq1$, the unique solution is $\ln\frac{8}{\delta}$.
Because by assumption~(\ref{eq:delta max}) $5<\ln\frac{8}{\delta}$,
hence $x_{0}>5$.

For $\sigma<1$, equation~(\ref{eq:v approx best condition log version again})
has at most one solution because the left-hand side is growing monotonically
with $x$ and is unbounded from above. However, because the highest
possible value of the left-hand side for $x=5$, which is realized
for $\sigma=\tau=0$, is smaller then the right hand side, $5-\frac{1}{2}\ln5.5<\ln\frac{8}{\delta}$,
there is exactly one solution $x_{0}>5$. 

We choose $v=\left\lceil x_{0}\right\rceil $, and determine $N$
using~(\ref{eq:N sufficient}). Because the solutions of~(\ref{eq:v approx best condition log version again})
increase monotonically with $x$, we have 
\begin{eqnarray}
\max\left(\frac{1-\sigma}{2},0\right)\ln\left(\frac{1}{2}+v+\tau\right)+\ln\frac{8}{\delta} & \leq & v\label{eq:v lower bound}\\
\max\left(\frac{1-\sigma}{2},0\right)\ln\left(\frac{1}{2}+v+\tau\right)+\ln\frac{8}{\delta}+1 & \geq & v\label{eq:v upper bound}
\end{eqnarray}
We now determine $E_{1}\left(s\right)$ and $D\left(s\right)$ to
accuracy $\frac{\delta}{3}$ and show that we can neglect within this
accuracy $E_{-1}\left(s\right)$, so that we can calculate $\zeta\left(s\right)$
to accuracy $\delta$.

\subsubsection{Upper bound for $E_{\mu}\left(m,s\right)$\label{sub:Upper bound for E_mu(m,s)}}

We express $E_{\mu}\left(m,s\right)$ as the remainder of a Taylor
expansion~(\cite{NIST_online},~1.4.35, 1.4.37), setting $f\left(z\right)=z^{s-1}$,
$a=m+\frac{i\mu}{2\pi N}$ and $b=m$, so that
\begin{eqnarray*}
E_{\mu}\left(m,s\right)=f\left(b\right)-\sum_{w=0}^{v-1}\frac{f^{\left(w\right)}\left(a\right)}{w!}\left(b-a\right)^{w} & = & \int_{a}^{b}\left(b-z\right)^{v-1}\frac{f^{\left(v\right)}\left(z\right)}{\left(v-1\right)!}dz
\end{eqnarray*}
Here we can assume that the integration runs over a straight line
segment from $a$ to $b$. The triangle inequality for integrals yields,
\begin{eqnarray*}
\left|\Gamma\left(1-s\right)E_{\mu}\left(m,s\right)\right| & \leq & \frac{\left(2\pi N\right)^{-v}}{\Gamma\left(v\right)}\left|\Gamma\left(1-s+v\right)\right|\max_{u\in\left[0;1\right]}\left|\left(m+i\frac{\mu u}{2\pi N}\right)^{s-1-v}\right|
\end{eqnarray*}
The last term is $\exp\left[\eta_{\mu}\left(m,s\right)\right]$ where
the real number $\eta_{\mu}\left(m,s\right)$ is at most
\begin{eqnarray*}
\eta_{\mu}\left(m,s\right) & = & \max_{u\in\left[0;1\right]}\Re\left[\left(\sigma-1-v+i\tau\right)\left(\ln\left|m+\frac{i\mu u}{2\pi N}\right|+i\mu\arctan\frac{u}{2\pi mN}\right)\right]\\
 & < & \left(\sigma-1-v\right)\ln m+\max_{u\in\left[0;1\right]}\left(-\mu\tau\arctan\frac{u}{2\pi mN}\right)\\
 & < & \left(\sigma-1-v\right)\ln m+\frac{1-\mu}{4}\pi\tau
\end{eqnarray*}
We have because of~(\ref{eq:0<=00003D sigma <=00003D 2}) and~(\ref{eq:v>=00003D5}),
\begin{eqnarray}
\left|1-s+v\right| & \geq & 4\label{eq:|1-s+v| >=00003D3}
\end{eqnarray}
Therefore, using the upper bounds~(\cite{NIST_online},~5.6.1 and~5.6.9),
\begin{eqnarray}
\Gamma\left(v\right)^{-1} & \leq & \left(2\pi\right)^{-\frac{1}{2}}v^{\frac{1}{2}-v}e^{v}\label{eq:Gamma(v) lower estimate NIST 5.6.1}\\
\left|\Gamma\left(1-s+v\right)\right| & < & e^{\frac{1}{6\left|1-s+v\right|}}\sqrt{2\pi}\left|\frac{1}{2}+v+\tau\right|^{\frac{1}{2}-\sigma+v}e^{-\frac{\pi}{2}\tau}\label{eq:Gamma function upper estimate NIST 5.6.9}
\end{eqnarray}
and the triangle inequality, we arrive with~(\ref{eq:|1-s+v| >=00003D3})
at the upper bound 
\begin{align}
 & \left|\left(2\pi\right)^{s-1}\Gamma\left(1-s\right)e^{i\mu\frac{\pi}{2}\left(1-s\right)}E_{\mu}\left(m,s\right)\right|\label{eq:E_mu(m,s) upper bound}\\
 & <e^{\frac{1}{24}}\left(2\pi\right)^{\sigma-1-v}N^{-v}v^{\frac{1}{2}-v}e^{v}\left(\frac{1}{2}+v+\tau\right)^{\frac{1}{2}-\sigma+v}\exp\left[\frac{\mu-1}{4}\pi\tau\right]m^{\sigma-v-1}\nonumber 
\end{align}

\subsubsection{Upper bound for $E_{-1}\left(s\right)$\label{sub:Upper bound for E_-1(s)}}

(\ref{eq:N sufficient}) implies $N>1$ and therefore 
\begin{eqnarray}
M & \geq & 2\label{eq:M>=00003D2}
\end{eqnarray}
For $\mu=-1$ and because of
\begin{eqnarray*}
\sum_{m=1}^{\infty}m^{\sigma-v-1} & \leq & \zeta\left(4\right)
\end{eqnarray*}
we can sum over all $m$ and have from~(\ref{eq:E_mu(m,s) upper bound})
and~(\ref{eq:M>=00003D2}) the upper bound,
\begin{eqnarray}
\left|E_{-1}\left(s\right)\right| & < & e^{\frac{1}{24}}\frac{\pi^{4}}{90}\left(2\pi\right)^{\sigma-1}\left(\frac{e}{4\pi}\right)^{v}v^{\frac{1}{2}-v}e^{\left[\left(v+\frac{1}{2}-\sigma\right)\ln\left(\frac{1}{2}+v+\tau\right)-\frac{\pi}{2}\tau\right]}\label{eq:E_-1 upper limit}
\end{eqnarray}
Because the term in square brackets decreases monotonically with $\tau$,
we choose its maximum value at $\tau=0$,
\begin{eqnarray*}
\left|E_{-1}\left(s\right)\right| & < & e^{\frac{1}{24}}\frac{\pi^{4}}{90}\left(2\pi\right)^{\sigma-1}\left(\frac{e}{4\pi}\right)^{v}\left(v+\frac{1}{2}\right)^{\frac{1}{2}}\left(v+\frac{1}{2}\right)^{\frac{1}{2}-\sigma}\left(1+\frac{1}{2v}\right)^{v}
\end{eqnarray*}
The last term is less than $\sqrt{e}$, so that
\begin{eqnarray*}
3\left|E_{-1}\left(s\right)\right| & < & \left[e^{\frac{1}{2}+\frac{1}{24}}\frac{\pi^{3}}{480}\left(\frac{v+\frac{1}{2}}{2\pi}\right)^{-\sigma}\left(v+\frac{1}{2}\right)\left(\frac{4\pi}{e^{2}}\right)^{-v}\right]\left(8e^{-v}\right)
\end{eqnarray*}
The term in square brackets is always less than one, because for $v=5$
it has because of $5.5<2\pi$ its maximum $<1$ for $\sigma=2$, and
for $v\geq6$ we have because of $6.5>2\pi$ an upper bound $<1$
by setting $\sigma=0$ and $v=6$. 

Hence we can neglect $E_{-1}\left(s\right)$ up to accuracy $\frac{\delta}{3}$,
because we see from~(\ref{eq:v lower bound}) that
\begin{eqnarray}
\ln\frac{8}{\delta} & \leq & v\label{eq:ln(8/delta)<v}
\end{eqnarray}
However, because for real argument $s$, $E_{1}\left(s\right)$ and
$E_{-1}\left(s\right)$ are complex conjugates, this shows that we
can neglect in this case both.

\subsubsection{Accuracy of $E_{1}\left(s\right)$\label{sub:Accuracy of  E_1(s)}}

In~(\ref{eq:Zetafast approximation}), we cut off the series for
$E_{1}\left(s\right)$ at $m=M$. We estimate now the rest $r_{E}$.
Using 
\begin{eqnarray*}
\sum_{m=M+1}^{\infty}m^{\sigma-1-v} & < & \int_{M}^{\infty}x^{\sigma-1-v}dx<\frac{M^{\sigma-v}}{v-2}<\frac{N^{\sigma-v}}{v-2}
\end{eqnarray*}
and the upper bound~(\ref{eq:E_mu(m,s) upper bound}) we have 
\begin{align}
\left|r_{E}\right|< & \frac{e^{\frac{1}{24}}}{8}\left(2\pi\right)^{\sigma-1}\frac{v^{\frac{1-\sigma}{2}}}{v-2}\left(\frac{e^{2}}{2\pi}\right)^{v}\left(1+\frac{\frac{1}{2}+\tau}{v}\right)^{v-\frac{\sigma}{2}}N^{\sigma-2v}\left[\left(\frac{1}{2}+v+\tau\right)^{\frac{1-\sigma}{2}}e^{-v}\frac{8}{\delta}\right]\delta\label{eq:|r_E| upper bound}
\end{align}
Because of $\frac{1-\sigma}{2}\leq\max\left(\frac{1-\sigma}{2},0\right)$
and~(\ref{eq:v lower bound}), we can replace the term in square
brackets by its upper bound $1$. Inserting $N$ of~(\ref{eq:N sufficient}),
we have because of
\begin{eqnarray*}
1.11^{2v-\sigma}> & \left(\frac{e^{2}}{2\pi}\right)^{\frac{5}{8}\left(2v-\sigma\right)}\geq & \left(\frac{e^{2}}{2\pi}\right)^{v}
\end{eqnarray*}
the upper bound
\begin{align}
\left|r_{E}\right|< & \left[\frac{3}{8}\frac{e^{\frac{1}{24}}}{v-2}\left(\frac{v}{4\pi^{2}}\right)^{\frac{1-\sigma}{2}}\right]\frac{\delta}{3}\label{eq:| r_E | bound}
\end{align}
The term in square brackets has for $v\leq4\pi^{2}$ its maximum at
$\sigma=2$ and $v=5$ as well as for $v\geq4\pi^{2}$ for $\sigma=0$
and $v=40$. In both cases this is less than one. Hence the value~(\ref{eq:N sufficient})
for $N$ is sufficient for the desired accuracy $\left|r_{E}\right|<\frac{\delta}{3}$.

\subsubsection{Accuracy of $D\left(s\right)$\label{sub:Accuracy of D(s)}}

At first we determine an upper bound for $N^{1-\sigma}$, using~(\ref{eq:N sufficient}),
\begin{eqnarray*}
\ln\frac{\left(\frac{N}{1.11}\right)^{1-\sigma}}{\delta}= & \left(-\ln8+\frac{\sigma-1}{2}\ln v\right) & +\left(\ln\frac{8}{\delta}+\frac{1-\sigma}{2}\ln\left(\frac{1}{2}+v+\tau\right)\right)
\end{eqnarray*}
The first term on the right hand side is at most
\begin{eqnarray*}
-\ln8+\frac{1}{2}\ln v & < & 0.003\cdot v
\end{eqnarray*}
The second term is because of~(\ref{eq:v lower bound}) and $\frac{1-\sigma}{2}\leq\max\left(\frac{1-\sigma}{2},0\right)$
at most $v$. Hence we have the upper bound
\begin{eqnarray}
N^{1-\sigma} & < & 1.11\,\delta e^{1.003v}\label{eq:N upper bound}
\end{eqnarray}
Using this bound and the triangle inequality, we get an upper bound
for $\left|r_{D}\right|$, 
\begin{eqnarray}
D\left(s\right) & = & \sum_{n=1}^{\left\lceil \lambda vN\right\rceil }n^{-s}Q\left(v,\frac{n}{N}\right)+r_{D}\label{eq:D(s) with rest term}\\
\left|r_{D}\right| & \leq & 1.11\,\delta e^{1.003v}\Gamma\left(v\right)^{-1}\sum_{n=\left\lceil \lambda vN\right\rceil +1}^{\infty}\left(\frac{n}{N}\right)^{-\sigma}\frac{1}{N}\int_{\frac{n}{N}}^{\infty}t^{v-1}e^{-t}dt\nonumber 
\end{eqnarray}
Because of $\sigma\geq0$, the integral is a decreasing function of
$\frac{n}{N}$, so that using inequality~(\ref{eq:Gamma(v) lower estimate NIST 5.6.1}),
we can bound the sum by the double integral
\begin{eqnarray}
\left|r_{D}\right| & \leq & 1.11\,\delta e^{1.003v}\left[\left(2\pi\right)^{-\frac{1}{2}}v^{\frac{1}{2}-v}e^{v}\right]\int_{\lambda v}^{\infty}du\int_{u}^{\infty}t^{v-1}e^{-t}dt\label{eq:r_D| prelim}
\end{eqnarray}
We have for $\lambda>1$ and $t\geq u\geq\lambda v$,
\begin{eqnarray*}
1 & < & \frac{1-\frac{v-1}{t}}{1-\frac{1}{\lambda}}
\end{eqnarray*}
and therefore 
\begin{eqnarray*}
\int_{u}^{\infty}t^{v-1}e^{-t}dt & < & \frac{\lambda}{\lambda-1}\int_{u}^{\infty}\left(t^{v-1}-\left(v-1\right)t^{v-2}\right)e^{-t}dt=\frac{\lambda}{\lambda-1}u^{v-1}e^{-u}
\end{eqnarray*}
and by the same argument
\begin{eqnarray*}
\int_{\lambda v}^{\infty}du\int_{u}^{\infty}t^{v-1}e^{-t}dt<\frac{\lambda}{\lambda-1}\int_{\lambda v}^{\infty}duu^{v-1}e^{-u} & < & \frac{\lambda}{\left(\lambda-1\right)^{2}}\lambda^{v}v^{v-1}e^{-\lambda v}
\end{eqnarray*}
Hence inequality~(\ref{eq:r_D| prelim}) becomes
\begin{eqnarray}
\left|r_{D}\right| & < & \frac{\delta}{3}\left[3\frac{\lambda}{\left(\lambda-1\right)^{2}}\frac{1.11}{\sqrt{2\pi v}}\right]e^{-v\left(\lambda-2.003-\ln\lambda\right)}\label{eq:| r_D | bound}
\end{eqnarray}
Assuming $\lambda\geq3$, the term in square brackets is always smaller
than one. Hence it suffices for an accuracy $\frac{\delta}{3}$ to
choose $\lambda=3.151$ because then 
\begin{eqnarray}
\lambda-2.003-\ln\lambda & > & 0\label{eq:lambda sufficient condition}
\end{eqnarray}

\subsubsection{Estimating the total number of summands\label{sub:Estimating the total number of summands}}

We have $\widehat{\lambda vN}<\lambda vN+1$ summands for $D\left(N,s\right)$
and $\left(v+1\right)M<\left(v+1\right)\left(N+1\right)$ summands
for $E_{1}\left(M,s\right)$. Therefore we have for the total number
of summands the upper bound 
\begin{eqnarray}
S & < & \left(\lambda v+v+1\right)N+v+2\label{eq:N_0 total number of summands}
\end{eqnarray}
Because of $v\geq5$ we have $1.11\left(\lambda v+v+1\right)<4.83v$
and hence from~(\ref{eq:N sufficient}) the upper bound
\begin{eqnarray}
S & < & 4.83\left[v\left(\frac{1}{2}+v+\tau\right)\right]^{\frac{1}{2}}+v+2\label{eq:S prelim}
\end{eqnarray}
We have for all $\tau\geq0$
\begin{eqnarray}
\tau-\frac{1}{2}\ln\left(2\tau\right) & > & \frac{3}{5}\tau
\end{eqnarray}
and therefore because of~(\ref{eq:tau large condition}) for all
$\sigma\geq0$, 
\begin{eqnarray*}
\tau-\frac{3}{2}-\max\left(\frac{1-\sigma}{2},0\right)\ln\left(\frac{1}{2}+\tau-\frac{3}{2}+\tau\right) & > & \ln\frac{8}{\delta}
\end{eqnarray*}
Because the solutions of~(\ref{eq:v approx best condition log version again})
increase monotonically with $x$, it follows 
\begin{eqnarray}
\tau & > & \frac{3}{2}+x_{0}>\frac{1}{2}+v\label{eq:t> v+1/2}
\end{eqnarray}
thus also $v<\left(v\tau\right)^{\frac{1}{2}}$, and therefore from~(\ref{eq:S prelim}),
\begin{eqnarray*}
S & < & \left(4.83\sqrt{2}+1\right)\sqrt{v}\sqrt{\tau}+2
\end{eqnarray*}
Using~(\ref{eq:v upper bound}) and~(\ref{eq:t> v+1/2}) we arrive
at the upper bound~(\ref{eq:Zetafast speed}).

\subsection{Proof of Theorem~\ref{thm:Zetafast algorithm for Dirichlet L-functions}\label{sub:Proof of Theorem=0000A0=00005Bthm:Zetafast Dirichlet=00005D}}

A Dirichlet $L$-function with primitive character $\chi\mod q$ is
given for $\sigma>0$ by
\begin{eqnarray}
L\left(s,\chi\right) & = & \sum_{n=1}^{\infty}\chi\left(n\right)n^{-s}\label{eq:Dirichlet L function}
\end{eqnarray}
and fulfills the functional equation~(\cite{NIST_online},~25.15.5
and~25.15.6),

\begin{align}
L\left(s,\chi\right) & =G\left(\chi\right)q^{-s}\left(2\pi\right)^{s-1}\Gamma\left(1-s\right)\sum_{\mu=\pm1}L\left(1-s,\overline{\chi}\right)\chi\left(-\mu\right)e^{i\mu\frac{\pi}{2}\left(1-s\right)}\label{eq:Dirichlet L-function functional equation}
\end{align}
where
\begin{eqnarray}
G\left(\chi\right) & = & \sum_{p=1}^{q}\chi\left(p\right)e^{\frac{2\pi pi}{q}}\label{eq:Gauss sum for chi}
\end{eqnarray}
is the Gauss sum. Repeating the arguments of section~\ref{sub:Proof of Theorem=0000A0=00005Bthm:Zetafast=00005D}
and assuming that $\chi$ is not the principal character, we have
at once its Zetafast algorithm~(\ref{eq:Zetafast Dirichlet L-function}).

\section*{Acknowledgment}

I would like to thank Henri Cohen for pointing out that the error
estimate could be strengthened.

\end{document}